\documentclass[12pt]{article}

\usepackage{amsmath,amssymb,amsthm,amscd,a4wide,upref}

\usepackage[enableskew]{youngtab}

\usepackage{hyperref}
\hypersetup{colorlinks,citecolor=blue,filecolor=black,linkcolor=blue,urlcolor=blue}
\usepackage{enumerate}
\usepackage{mathtools}
\usepackage{enumitem}

\usepackage{lmodern}     
\usepackage[T1]{fontenc}

\makeatletter
\DeclareTextCommand{\textprime}{\encodingdefault}{%
  \mbox{$\m@th'\kern-\scriptspace$}%
}
\makeatother

\begin{document}


\newcommand{\ad}{{\rm ad}}
\newcommand{\cri}{{\rm cri}}
\newcommand{\ext}{{\rm ext}}
\newcommand{\row}{{\rm row}}
\newcommand{\col}{{\rm col}}
\newcommand{\End}{{\rm{End}\ts}}
\newcommand{\Rep}{{\rm{Rep}\ts}}
\newcommand{\Hom}{{\rm{Hom}}}
\newcommand{\Mat}{{\rm{Mat}}}
\newcommand{\ch}{{\rm{ch}\ts}}
\newcommand{\chara}{{\rm{char}\ts}}
\newcommand{\diag}{{\rm diag}}
\newcommand{\st}{{\rm st}}
\newcommand{\non}{\nonumber}
\newcommand{\wt}{\widetilde}
\newcommand{\wh}{\widehat}
\newcommand{\ol}{\overline}
\newcommand{\ot}{\otimes}
\newcommand{\la}{\lambda}
\newcommand{\La}{\Lambda}
\newcommand{\De}{\Delta}
\newcommand{\al}{\alpha}
\newcommand{\be}{\beta}
\newcommand{\ga}{\gamma}
\newcommand{\Ga}{\Gamma}
\newcommand{\ep}{\epsilon}
\newcommand{\ka}{\kappa}
\newcommand{\vk}{\varkappa}
\newcommand{\vt}{\vartheta}
\newcommand{\si}{\sigma}
\newcommand{\vs}{\varsigma}
\newcommand{\vp}{\varphi}
\newcommand{\de}{\delta}
\newcommand{\ze}{\zeta}
\newcommand{\om}{\omega}
\newcommand{\Om}{\Omega}
\newcommand{\ee}{\epsilon^{}}
\newcommand{\su}{s^{}}
\newcommand{\hra}{\hookrightarrow}
\newcommand{\ve}{\varepsilon}
\newcommand{\ts}{\,}
\newcommand{\vac}{\mathbf{1}}
\newcommand{\vacu}{|0\rangle}
\newcommand{\di}{\partial}
\newcommand{\qin}{q^{-1}}
\newcommand{\tss}{\hspace{1pt}}
\newcommand{\Sr}{ {\rm S}}
\newcommand{\U}{ {\rm U}}
\newcommand{\BL}{ {\overline L}}
\newcommand{\BE}{ {\overline E}}
\newcommand{\BP}{ {\overline P}}
\newcommand{\AAb}{\mathbb{A}\tss}
\newcommand{\CC}{\mathbb{C}\tss}
\newcommand{\KK}{\mathbb{K}\tss}
\newcommand{\QQ}{\mathbb{Q}\tss}
\newcommand{\SSb}{\mathbb{S}\tss}
\newcommand{\TT}{\mathbb{T}\tss}
\newcommand{\ZZ}{\mathbb{Z}\tss}
\newcommand{\DY}{ {\rm DY}}
\newcommand{\X}{ {\rm X}}
\newcommand{\Y}{ {\rm Y}}
\newcommand{\Z}{{\rm Z}}
\newcommand{\Ac}{\mathcal{A}}
\newcommand{\Lc}{\mathcal{L}}
\newcommand{\Mc}{\mathcal{M}}
\newcommand{\Pc}{\mathcal{P}}
\newcommand{\Qc}{\mathcal{Q}}
\newcommand{\Rc}{\mathcal{R}}
\newcommand{\Sc}{\mathcal{S}}
\newcommand{\Tc}{\mathcal{T}}
\newcommand{\Bc}{\mathcal{B}}
\newcommand{\Cc}{\mathcal{C}}
\newcommand{\Ec}{\mathcal{E}}
\newcommand{\Fc}{\mathcal{F}}
\newcommand{\Gc}{\mathcal{G}}
\newcommand{\Hc}{\mathcal{H}}
\newcommand{\Uc}{\mathcal{U}}
\newcommand{\Vc}{\mathcal{V}}
\newcommand{\Wc}{\mathcal{W}}
\newcommand{\Xc}{\mathcal{X}}
\newcommand{\Yc}{\mathcal{Y}}
\newcommand{\Ar}{{\rm A}}
\newcommand{\Br}{{\rm B}}
\newcommand{\Ir}{{\rm I}}
\newcommand{\Fr}{{\rm F}}
\newcommand{\Jr}{{\rm J}}
\newcommand{\Or}{{\rm O}}
\newcommand{\GL}{{\rm GL}}
\newcommand{\Spr}{{\rm Sp}}
\newcommand{\Rr}{{\rm R}}
\newcommand{\Zr}{{\rm Z}}
\newcommand{\gl}{\mathfrak{gl}}
\newcommand{\middd}{{\rm mid}}
\newcommand{\ev}{{\rm ev}}
\newcommand{\Pf}{{\rm Pf}}
\newcommand{\Norm}{{\rm Norm\tss}}
\newcommand{\oa}{\mathfrak{o}}
\newcommand{\spa}{\mathfrak{sp}}
\newcommand{\osp}{\mathfrak{osp}}
\newcommand{\f}{\mathfrak{f}}
\newcommand{\g}{\mathfrak{g}}
\newcommand{\h}{\mathfrak h}
\newcommand{\n}{\mathfrak n}
\newcommand{\m}{\mathfrak m}
\newcommand{\z}{\mathfrak{z}}
\newcommand{\Zgot}{\mathfrak{Z}}
\newcommand{\p}{\mathfrak{p}}
\newcommand{\sll}{\mathfrak{sl}}
\newcommand{\whg}{\wh{\g}}
\newcommand{\gll}{\g^{}_{\ell}}
\newcommand{\gllh}{\wh{\g}^{}_{\ell}}
\newcommand{\gllm}{\g^{}_{\ell,\ell'}}
\newcommand{\glls}{\g^*_{\ell}}
\newcommand{\agot}{\mathfrak{a}}
\newcommand{\bgot}{\mathfrak{b}}
\newcommand{\qdet}{ {\rm qdet}\ts}
\newcommand{\tra}{ {\rm t}}
\newcommand{\Ber}{ {\rm Ber}\ts}
\newcommand{\HC}{ {\mathcal HC}}
\newcommand{\cdet}{{\rm cdet}}
\newcommand{\rdet}{{\rm rdet}}
\newcommand{\tr}{ {\rm tr}}
\newcommand{\gr}{ {\rm gr}\ts}
\newcommand{\str}{ {\rm str}}
\newcommand{\loc}{{\rm loc}}
\newcommand{\Fun}{{\rm{Fun}\ts}}
\newcommand{\Gr}{{\rm G}}
\newcommand{\sgn}{ {\rm sgn}\ts}
\newcommand{\sign}{{\rm sgn}}
\newcommand{\ba}{\bar{a}}
\newcommand{\bb}{\bar{b}}
\newcommand{\bi}{\bar{\imath}}
\newcommand{\bj}{\bar{\jmath}}
\newcommand{\bk}{\bar{k}}
\newcommand{\bl}{\bar{l}}
\newcommand{\hb}{\mathbf{h}}
\newcommand{\Sym}{\mathfrak S}
\newcommand{\fand}{\quad\text{and}\quad}
\newcommand{\Fand}{\qquad\text{and}\qquad}
\newcommand{\For}{\qquad\text{or}\qquad}
\newcommand{\for}{\quad\text{or}\quad}
\newcommand{\grpr}{{\rm gr}^{\tss\prime}\ts}
\newcommand{\degpr}{{\rm deg}^{\tss\prime}\tss}
\newcommand{\bideg}{{\rm bideg}\ts}

\renewcommand{\theequation}{\arabic{section}.\arabic{equation}}

\numberwithin{equation}{section}

\newtheorem{thm}{Theorem}[section]
\newtheorem{lem}[thm]{Lemma}
\newtheorem{prop}[thm]{Proposition}
\newtheorem{cor}[thm]{Corollary}
\newtheorem{conj}[thm]{Conjecture}
\newtheorem*{mthm}{Main Theorem}
\newtheorem*{mthma}{Theorem A}
\newtheorem*{mthmb}{Theorem B}
\newtheorem*{mthmc}{Theorem C}
\newtheorem*{mthmd}{Theorem D}

\theoremstyle{definition}
\newtheorem{defin}[thm]{Definition}

\theoremstyle{remark}
\newtheorem{remark}[thm]{Remark}
\newtheorem{example}[thm]{Example}
\newtheorem{examples}[thm]{Examples}

\newcommand{\bth}{\begin{thm}}
\renewcommand{\eth}{\end{thm}}
\newcommand{\bpr}{\begin{prop}}
\newcommand{\epr}{\end{prop}}
\newcommand{\ble}{\begin{lem}}
\newcommand{\ele}{\end{lem}}
\newcommand{\bco}{\begin{cor}}
\newcommand{\eco}{\end{cor}}
\newcommand{\bde}{\begin{defin}}
\newcommand{\ede}{\end{defin}}
\newcommand{\bex}{\begin{example}}
\newcommand{\eex}{\end{example}}
\newcommand{\bes}{\begin{examples}}
\newcommand{\ees}{\end{examples}}
\newcommand{\bre}{\begin{remark}}
\newcommand{\ere}{\end{remark}}
\newcommand{\bcj}{\begin{conj}}
\newcommand{\ecj}{\end{conj}}

\newcommand{\bal}{\begin{aligned}}
\newcommand{\eal}{\end{aligned}}
\newcommand{\beq}{\begin{equation}}
\newcommand{\eeq}{\end{equation}}
\newcommand{\ben}{\begin{equation*}}
\newcommand{\een}{\end{equation*}}

\newcommand{\bpf}{\begin{proof}}
\newcommand{\epf}{\end{proof}}

\def\beql#1{\begin{equation}\label{#1}}

\newcommand{\Res}{\mathop{\mathrm{Res}}}

\title{\Large\bf Eigenvalues of quantum Gelfand invariants}

\author{{Naihuan Jing,\quad Ming Liu\quad and\quad Alexander Molev}}

\date{} 
\maketitle


\begin{abstract}
We consider the quantum Gelfand invariants which first appeared
in a landmark paper
by Reshetikhin, Takhtadzhyan and Faddeev (1989).
We calculate the eigenvalues of the invariants acting in irreducible highest weight
representations of the quantized enveloping algebra for $\gl_n$.
The calculation is based on Liouville-type formulas relating
two families of
central elements in the quantum affine algebras of type $A$.



\end{abstract}

\section{Introduction}
\label{sec:int}

The {\em quantized enveloping algebras} and {\em quantum affine algebras} associated with simple
Lie algebras comprise remarkable families of
quantum groups, as introduced by V.~Drinfeld~\cite{d:ha}
and M.~Jimbo~\cite{j:qd}. These algebras and their representations
have since found numerous connections with
many areas in mathematics and physics.

In this paper we will be concerned with those families associated with the general linear
Lie algebras $\gl_n$. Both the quantized enveloping algebra $\U_q(\gl_n)$ and
the quantum affine algebra $\U_q(\wh\gl_n)$ admit $R$-{\em matrix} (or $RTT$) {\em presentations}
going back to the work of the Leningrad school headed by L.~D.~Faddeev; see e.g \cite{ks:qs}
and \cite{rtf:ql}
for reviews of the foundations of the $R$-matrix approach
originated in the quantum inverse scattering method.

Central elements in both $\U_q(\gl_n)$ and $\U_q(\wh\gl_n)$ are constructed with the use
of the $R$-matrix presentations and found as coefficients of
the respective {\em quantum determinants}; see
\cite{c:ni}, \cite{j:qu}, \cite{ks:qs}, and also \cite{e:ce}
for more general constructions of central elements in the quantized enveloping algebras and
quantum affine algebras.
As pointed out in \cite{rtf:ql}, the quantum traces of powers of generator matrices are central
in $\U_q(\gl_n)$; see also \cite{b:cq}. By taking the limit $q\to 1$ one recovers
the central
elements of $\U(\gl_n)$ going back to \cite{g:tc}, which
are known as the {\em Gelfand invariants}.
Note that a different
generalization of the Gelfand invariants for $\gl_n$
as central elements in $\U_q(\gl_n)$,
was given in \cite{gzb:gg}, where their eigenvalues in irreducible highest weight
representations were calculated.

A new family of central elements in $\U_q(\wh\gl_n)$ was given in \cite{br:nb}
and they were
related to the quantum determinants by Liouville-type formulas, although
they were not accompanied by proofs. This result is quite
analogous to the corresponding quantum Liouville formulas for the Yangians
originated in \cite{n:qb} and we give a complete proof in this paper.

By taking the images of the new central elements of $\U_q(\wh\gl_n)$ under the evaluation homomorphism,
we recover the quantum Gelfand invariants of \cite{rtf:ql}. The Liouville formulas will then
allow us to calculate the eigenvalues of the quantum Gelfand invariants
in irreducible highest weight representations of $\U_q(\gl_n)$. We thus obtain $q$-analogues
of the Perelomov--Popov formulas \cite{pp:co}. To recall the eigenvalue formulas from \cite{pp:co},
consider the irreducible highest weight representation $L(\la)$ of $\gl_n$ with the highest weight
$\la=(\la_1,\dots,\la_n)$ and combine the standard basis elements $E_{ij}$ into the
matrix $E=[E_{ij}]$.
Then the eigenvalue of the Gelfand invariant $\tr\ts E^m$ in $L(\la)$ is found by
\beql{ppexp}
\tr\ts E^m\mapsto\sum_{k=1}^n\ts \ell^{\ts m}_k\ts \frac{(\ell_1-\ell_k+1)\dots (\ell_n-\ell_k+1)}
{(\ell_1-\ell_k)\ldots\wedge\dots (\ell_n-\ell_k)},
\eeq
where $\ell_i=\la_i+n-i$ and the symbol $\wedge$ indicates that the zero factor is skipped.
Formula \eqref{ppexp} can be derived with the use of $R$-matrix calculations
in the Yangian $\Y(\gl_n)$ or in the universal enveloping algebra $\U(\gl_n)$;
see \cite[Sec.~7.1]{m:yc} and \cite[Sec.~4.8]{m:so}, respectively.

Our main result concerning the quantum Gelfand invariants
is the following theorem, where we calculate the
eigenvalues of the quantum traces $\tr_q\tss M^m$ of the powers of the generator
matrix $M=L^{-}(L^{+})^{-1}$ in the
representation $L_q(\la)$
of $\U_q(\gl_n)$ (see Section~\ref{sec:qgi} for the definitions). We use a standard notation
for the $q$-numbers
\ben
[k]_q=\frac{q^k-q^{-k}}{q-\qin},\qquad k\in\ZZ.
\een

\bth\label{thm:eigen}
The eigenvalue of the quantum Gelfand invariant $\tr_q\tss M^m$ in $L_q(\la)$ is found by
\beql{qppexp}
\tr_q\tss M^m\mapsto\sum_{k=1}^n\ts q^{2 \ell_k m}\ts \frac{[\ell_1-\ell_k+1]_q\dots [\ell_n-\ell_k+1]_q}
{[\ell_1-\ell_k]_q\ldots\wedge\dots [\ell_n-\ell_k]_q}.
\eeq
\eth

We will prove Theorem~\ref{thm:eigen} in Section~\ref{sec:qgi} by deriving it from the
Liouville formula given in Theorem~\ref{thm:liov} in a way similar to \cite[Sec.~7.1]{m:yc}.
We also consider three more families of central elements of $\U_q(\gl_n)$
which, however, turn out to coincide with the quantum Gelfand invariants, up to a possible
replacement $q\mapsto\qin$.

In the limit $q\to 1$ we have
\beql{cllim}
\frac{M-1}{q-\qin}\rightarrow E
\eeq
by \eqref{taugen1} and \eqref{taugen2} below. Therefore,
the Gelfand invariants in $\U(\gl_n)$
are recovered from the elements $\tr_q\tss M^m$ in the limit, while
the Perelomov--Popov formulas \eqref{ppexp} follow from \eqref{qppexp}; see
Remark~\ref{rem:pp} below.

We point out a related recent work \cite{jw:cr}, where
explicit formulas for certain central elements in the reflection equation algebras
were given and their relation with
the quantum Gelfand invariants were reviewed. This includes
the connection with an earlier construction of central elements in
\cite{r:qh} and with
the Cayley--Hamilton theorem and Newton identities
of \cite{gps:hs}, \cite{iop:qm} and \cite{nt:yg}.

\medskip

This work was completed during the first and third named authors' visits to the South China University
of Technology and to the Shanghai University. They are grateful to the Departments of Mathematics
in both universities for the warm hospitality.

\section{Liouville formulas}
\label{sec:lf}

We will regard $q$ as a nonzero complex number which is not
a root of unity. Recall the
$R$-matrix presentation of the quantum affine algebra $\U_q(\wh\gl_n)$ as introduced in \cite{rs:ce}.
We follow \cite{df:it} and use the same settings as in our earlier work \cite{jlm:qs}.
Let $e_{ij}\in \End\CC^n$ denote the standard matrix units. Consider the $R$-matrix
\beql{rf}
R(x)=\frac{f(x)}{q-\qin x}(R-x\tss \wt R),
\eeq
where
\beql{R}
R=q\ts\sum_{i}e_{ii}\ot e_{ii}+\sum_{i\ne j}e_{ii}\ot e_{jj}
+(q-\qin)\sum_{i< j}e_{ij}\ot e_{ji}
\eeq
and
\ben
\wt R=\qin\sum_{i}e_{ii}\ot e_{ii}+\sum_{i\ne j}e_{ii}\ot e_{jj}
-(q-\qin)\sum_{i> j}e_{ij}\ot e_{ji},
\een
while the formal power series
\ben
f(x)=1+\sum_{k=1}^{\infty}f_kx^k,\qquad f_k=f_k(q),
\een
is uniquely determined
by the relation
\ben
f(xq^{2n})=f(x)\ts\frac{(1-xq^2)\tss(1-xq^{2n-2})}{(1-x)\tss(1-xq^{2n})}.
\een

The {\em quantum affine algebra $\U_q(\wh\gl_n)$}
is generated by elements
\ben
l^+_{ij}[-r],\qquad l^-_{ij}[r]\qquad\text{with}\quad 1\leqslant i,j\leqslant n,\qquad r=0,1,\dots,
\een
and the invertible central element $q^c$,
subject to the defining relations
\begin{align}
l^+_{ji}[0]&=l^-_{ij}[0]=0\qquad&&\text{for}\qquad 1\leqslant i<j\leqslant n,
\non\\
l^+_{ii}[0]\ts l^-_{ii}[0]&=l^-_{ii}[0]\ts l^+_{ii}[0]=1\qquad&&\text{for}\qquad i=1,\dots,n,
\non
\end{align}
and
\begin{align}
R(u/v)L_1^{\pm}(u)L_2^{\pm}(v)&=L_2^{\pm}(v)L_1^{\pm}(u)R(u/v),
\label{RLL}\\[0.2em]
R(uq^{-c}/v)L_1^{+}(u)L_2^{-}(v)&=L_2^{-}(v)L_1^{+}(u)R(uq^{c}/v).
\label{RLLpm}
\end{align}
In the last two formulas we consider the matrices $L^{\pm}(u)=\big[\tss l^{\pm}_{ij}(u)\big]$,
whose entries are formal power series in $u$ and $u^{-1}$,
\ben
l^{+}_{ij}(u)=\sum_{r=0}^{\infty}l^{+}_{ij}[-r]\tss u^r,\qquad
l^{-}_{ij}(u)=\sum_{r=0}^{\infty}l^{-}_{ij}[r]\tss u^{-r}.
\een
Here and below we regard the matrices as elements
\ben
L^{\pm}(u)=\sum_{i,j=1}^n e_{ij}\ot l^{\pm}_{ij}(u)\in\End\CC^n\ot\U_q(\wh\gl_n)[[u^{\pm1}]]
\een
and use a subscript to indicate a copy of the matrix in the multiple
tensor product algebra
\beql{multtpr}
\underbrace{\End\CC^n\ot\dots\ot\End\CC^n}_k\ot\U_q(\wh\gl_n)[[u^{\pm1}]]
\eeq
so that
\ben
L^{\pm}_a(u)=\sum_{i,j=1}^n 1^{\ot (a-1)}\ot e_{ij}\ot 1^{\ot (k-a)}\ot l^{\pm}_{ij}(u).
\een

We regard the usual matrix transposition also as the linear map
\ben
t:\End\CC^n\to\End\CC^n,\qquad e_{ij}\mapsto e_{ji}.
\een
For any $a\in\{1,\dots,k\}$ we will denote by $t_a$ the corresponding
partial transposition on the algebra \eqref{multtpr} which acts as $t$ on the
$a$-th copy of $\End \CC^n$ and as the identity map on all the other tensor factors.

The scalar factor in \eqref{rf} is necessary for the $R$-matrix to satisfy
the {\em crossing symmetry relations} \cite{fri:qa}. We will use one of them
given by
\beql{cs}
\big(R_{12}(x)^{-1}\big)^{t_2} D_2 R_{12}(xq^{2n})^{t_2}=D_2,
\eeq
where $D$ denotes the diagonal $n\times n$ matrix
\ben
D=\diag\big[q^{n-1}, q^{n-3},\dots, q^{-n+1}\big].
\een

The following proposition was stated in \cite[Eq. (4.28)]{br:nb} without proof.

\bpr\label{prop:equiz}
There exist a series $z^+(u)$ in $u$ and a series $z^-(u)$ in $u^{-1}$ with coefficients
in the algebra $\U_q(\wh\gl_n)$ such that
\begin{align}\label{lltr}
L^{\pm}(uq^{2n})^tD(L^{\pm}(u)^{-1})^t&=z^{\pm}(u)D \\
\intertext{and}
(L^{\pm}(u)^{-1})^tD^{-1} L^{\pm}(uq^{2n})^t&=z^{\pm}(u)D^{-1}.
\label{secllrt}
\end{align}
Moreover, the coefficients of the series $z^{\pm}(u)$ belong to the center of the algebra
$\U_q(\wh\gl_n)$.
\epr

\bpf
Multiply both sides of \eqref{RLL} by $L^{\pm}_2(v)^{-1}$ from the left and the right
and apply the transposition $t_2$ to get
\ben
R(u/v)^{t_2}\big(L^{\pm}_2(v)^{-1}\big)^tL^{\pm}_1(u)
=L^{\pm}_1(u)\big(L^{\pm}_2(v)^{-1}\big)^tR(u/v)^{t_2}
\een
and hence
\beql{rtrain}
\big(R(u/v)^{t_2}\big)^{-1}L^{\pm}_1(u)\big(L^{\pm}_2(v)^{-1}\big)^t=
\big(L^{\pm}_2(v)^{-1}\big)^tL^{\pm}_1(u)\big(R(u/v)^{t_2}\big)^{-1}.
\eeq
Use the crossing symmetry relation \eqref{cs} to replace the $R$-matrix by
\ben
\big(R(u/v)^{t_2}\big)^{-1}=D_2^{-1}\big(R(u/vq^{2n})^{-1}\big)^{t_2}D_2
\een
and get
\ben
\big(R(u/vq^{2n})^{-1}\big)^{t_2}D_2L^{\pm}_1(u)\big(L^{\pm}_2(v)^{-1}\big)^tD_2^{-1}=
D_2\big(L^{\pm}_2(v)^{-1}\big)^tL^{\pm}_1(u)D_2^{-1}\big(R(u/vq^{2n})^{-1}\big)^{t_2}.
\een
Now cancel the scalar factors appearing in \eqref{rf} on both sides of this relation and observe
that the $R$-matrix $R-x\tss \wt R$ evaluated at $x=1$ equals $(q-\qin)P$, where $P$
is the permutation operator. Therefore,
\ben
\big((R-\wt R)^{-1}\big)^{t_2}=\frac{1}{q-\qin}\ts Q\qquad\text{with}\quad
Q=\sum_{i,j=1}^n e_{ij}\ot e_{ij}.
\een
Hence, by taking $u=vq^{2n}$ we get
\ben
Q D_2L^{\pm}_1(vq^{2n})\big(L^{\pm}_2(v)^{-1}\big)^tD_2^{-1}
=D_2\big(L^{\pm}_2(v)^{-1}\big)^tL^{\pm}_1(vq^{2n})D_2^{-1}Q.
\een
Since $Q$ is an operator in $\End\CC^n\ot\End\CC^n$ with a one-dimensional image,
both sides must be equal to $Q\tss z^{\ts\pm}(v)$ for series $z^{\ts\pm}(v)$
with coefficients in the quantum affine algebra.
Using the relations $QX_1=QX_{2}^t$ and $X_1 Q=X_{2}^t Q$ which hold for an arbitrary matrix $X$, we can
write the definition of $z^{\ts\pm}(v)$ as
\begin{align}\label{qld}
QL^{\pm}_2(vq^{2n})^tD_2\big(L^{\pm}_2(v)^{-1}\big)^t&=Q\tss D_2\tss z^{\ts\pm}(v)\\
\intertext{and}
\big(L^{\pm}_2(v)^{-1}\big)^t D_2^{-1}L^{\pm}_2(vq^{2n})^tQ&=D_2^{-1}\tss Q\tss z^{\ts\pm}(v).
\label{opqld}
\end{align}
By taking trace over the first copy of $\End\CC^n$ on both sides of \eqref{qld}
and \eqref{opqld} we arrive at
\eqref{lltr} and \eqref{secllrt}, respectively.

We will now use \eqref{secllrt}
to show that the series $z^-(v)$ commutes with $L^+(u)$.
We have
\beql{loz}
L^+_1(u)z^-(v)=L^+_1(u)D_2\big(L^{-}_2(v)^{-1}\big)^t D_2^{-1}L^{-}_2(vq^{2n})^t.
\eeq
Transform relation \eqref{RLLpm} in the same way as we did for \eqref{RLL}
in the beginning of the proof to get the following counterpart of \eqref{rtrain}:
\ben
\big(R(uq^{-c}/v)^{t_2}\big)^{-1}L^{+}_1(u)\big(L^{-}_2(v)^{-1}\big)^t=
\big(L^{-}_2(v)^{-1}\big)^tL^{+}_1(u)\big(R(uq^{c}/v)^{t_2}\big)^{-1}.
\een
Hence, the right hand side of \eqref{loz} equals
\beql{drt}
D_2 R_{12}(uq^{-c}/v)^{t_2}\big(L^{-}_2(v)^{-1}\big)^t L^+_1(u)
\big(R_{12}(uq^{c}/v)^{t_2}\big)^{-1}D_2^{-1}L^{-}_2(vq^{2n})^t.
\eeq
Applying again \eqref{cs}, we can write
\ben
\big(R_{12}(uq^c/v)^{t_2}\big)^{-1}D_2^{-1}=D_2^{-1}\big(R_{12}(uq^c/vq^{2n})^{-1}\big)^{t_2}.
\een
Continue transforming \eqref{drt}
by using the following consequence of \eqref{RLLpm}:
\ben
L^+_1(u)\big(R_{12}(uq^c/vq^{2n})^{-1}\big)^{t_2}L^{-}_2(vq^{2n})^t
=L^{-}_2(vq^{2n})^t\big(R_{12}(uq^{-c}/vq^{2n})^{-1}\big)^{t_2}L^{-}_2(vq^{2n})^t,
\een
so that \eqref{drt} becomes
\ben
D_2 R_{12}(uq^{-c}/v)^{t_2}\big(L^{-}_2(v)^{-1}\big)^t D_2^{-1}L^{-}_2(vq^{2n})^t
\big(R_{12}(uq^{-c}/vq^{2n})^{-1}\big)^{t_2} L^+_1(u).
\een
By \eqref{secllrt} this simplifies to
\ben
D_2 R_{12}(uq^{-c}/v)^{t_2}D_2^{-1}
\big(R_{12}(uq^{-c}/vq^{2n})^{-1}\big)^{t_2}z^-(v) L^+_1(u)
\een
which equals $z^-(v)L^+_1(u)$ by \eqref{cs}.
This proves that $L^+_1(u)z^-(v)=z^-(v)L^+_1(u)$. The relation $L^-_1(u)z^-(v)=z^-(v)L^-_1(u)$
and the centrality of $z^+(v)$ are verified in the same way.
\epf

\bre\label{rem:bcd}
The counterparts of the series $z^{\pm}(u)$
for the quantum affine algebras of types $B$, $C$ and $D$
appear in \cite[Prop.~3.3]{jlm:ib-c}
and \cite[Prop.~3.3]{jlm:ib-bd}, where they were introduced
by
relations analogous to \eqref{lltr}
and \eqref{secllrt}.
\ere

\bco\label{cor:zu}
We have the formulas
\begin{align}\label{zu}
z^{\pm}(u)&=\frac{1}{[n]_q}\ts \tr\ts DL^{\pm}(uq^{2n})L^{\pm}(u)^{-1}\\
\intertext{and}
z^{\pm}(u)&=\frac{1}{[n]_q}\ts \tr\ts D^{-1}L^{\pm}(u)^{-1}L^{\pm}(uq^{2n}).
\label{opzu}
\end{align}
\eco

\bpf
The formulas follow by taking trace on both sides of the respective matrix relations
\eqref{lltr} and \eqref{secllrt}.
\epf

Recall that the {\em quantum determinants} $\qdet L^{+}(u)$ and $\qdet L^{-}(u)$
are series in $u$ and $u^{-1}$, respectively, whose coefficients belong to the center
of the quantum affine algebra $\U_q(\wh\gl_n)$:
\ben
\qdet L^{\pm}(u)=
\sum_{\si\in \Sym_n} (-q)^{-l(\si)} \ts l^{\pm}_{\si(1)1}(uq^{2n-2})\cdots
l^{\pm}_{\si(n)n}(u),
\een
where $l(\si)$ denotes the length of the permutation $\si$.

The following is a $q$-analogue of the quantum Liouville formula of \cite{n:qb}.
It was stated in \cite[Eq. (4.32)]{br:nb}.

\bth\label{thm:liov}
We have the relations
\beql{liouvpm}
z^{\pm}(u)=\frac{\qdet L^{\pm}(uq^{2})}{\qdet L^{\pm}(u)}.
\eeq
\eth

\bpf
We follow \cite[Sec.~1.9]{m:yc} and use the {\em quantum comatrices} ${\wh L}^{\pm}(u)$
introduced in \cite{nt:yg} and \cite{t:cm}.
They are defined by the relations
\beql{comatrix}
{\wh L}^{\pm}(uq^2)\tss L^{\pm}(u)=\qdet L^{\pm}(u)\ts 1,
\eeq
where $1$ denotes the identity matrix. We will derive formulas
for the entries of the matrices ${\wh L}^{\pm}(u)$ by using
the quantum minor relations reviewed e.g. in \cite{hm:qa}
which we outline below.

Recall that the $q$-permutation operator
$P^{\tss q}\in\End\CC^n\ot\End\CC^n$ is
defined by
\ben
P^{\tss q}=\sum_{i=1}^n e_{ii}\ot e_{ii}+ q\tss\sum_{i> j}e_{ij}\ot
e_{ji}+ \qin\sum_{i< j}e_{ij}\ot e_{ji}.
\een
The action of the symmetric group $\Sym_k$ on the space $(\CC^n)^{\ot\tss k}$
can be defined by setting $s_i\mapsto P^{\tss q}_{s_i}:=
P^{\tss q}_{i,i+1}$ for $i=1,\dots,k-1$,
where $s_i$ denotes the transposition $(i,i+1)$.
If $\si=s_{i_1}\cdots s_{i_l}$ is a reduced decomposition
of an element $\si\in \Sym_k$ then we set
$P^{\tss q}_{\si}=P^{\tss q}_{s_{i_1}}\cdots P^{\tss q}_{s_{i_l}}$.
Denote by $e_1,\dots,e_n$ the canonical
basis vectors of $\CC^n$. Then for any indices $a_1<\dots<a_k$
and any $\tau\in \Sym_k$ we have
\beql{pq}
P^{\tss q}_{\si}(e_{a_{\tau(1)}}\ot\dots\ot e_{a_{\tau(k)}})
=q^{l(\si\tau^{-1})-l(\tau)}\ts
e_{a_{\tau\si^{-1}(1)}}\ot\dots\ot e_{a_{\tau\si^{-1}(k)}}.
\eeq
We denote by $A^{(k)}_q$ the $q$-antisymmetrizer
\beql{antisym}
A^{(k)}_q=\sum_{\si\in\Sym_k}\sgn\ts\si\cdot P^{\tss q}_{\si}.
\eeq
The defining relations \eqref{RLL} and the {\em fusion procedure} \cite{c:ni} for
the quantum affine algebra imply the relations
\beql{anitt}
A^{(k)}_q L^{\pm}_1(uq^{2k-2})\cdots L^{\pm}_k(u)=L^{\pm}_k(u)\cdots L^{\pm}_1(uq^{2k-2})\tss  A^{(k)}_q.
\eeq
The quantum minors are the series ${l^{\pm}\ts}^{a_1\cdots\ts a_k}_{b_1\cdots\ts b_k}(u)$
with coefficients in $\U_q(\wh\gl_n)$
defined by the expansion
of the elements \eqref{anitt}
as
\ben
\sum_{a_i,b_i}
e_{a_1b_1}\ot\cdots\ot e_{a_kb_k}\ot {l^{\pm}\ts}^{a_1\cdots\ts a_k}_{b_1\cdots\ts b_k}(u).
\een
Explicit formulas for the quantum minors have the form: if $a_1<\cdots<a_k$ then
\beql{qminorgen}
{l^{\pm}\ts}^{a_1\cdots\ts a_k}_{b_1\cdots\ts b_k}(u)=
\sum_{\si\in \Sym_k} (-q)^{-l(\si)} \cdot l^{\pm}_{a_{\si(1)}b_1}(uq^{2k-2})\cdots
l^{\pm}_{a_{\si(k)}b_k}(u).
\eeq
If $b_1<\cdots<b_k$ (and the $a_i$ are arbitrary) then
\ben
{l^{\pm}\ts}^{a_1\cdots\ts a_k}_{b_1\cdots\ts b_k}(u)=
\sum_{\si\in \Sym_k} (-q)^{l(\si)} \cdot l^{\pm}_{a_kb_{\si(k)}}(u)\cdots
l^{\pm}_{a_1b_{\si(1)}}(uq^{2k-2})
\een
and for any $\tau\in\Sym_k$ we have
\ben
{l^{\pm}\ts}^{a_1\cdots\ts a_k}_{b_{\tau(1)}\cdots\ts b_{\tau(k)}}(u)=
(-q)^{-l(\tau)}{l^{\pm}\ts}^{a_1\cdots\ts a_k}_{b_1\cdots\ts b_k}(u).
\een

The following lemma was pointed out in \cite{nt:yg} and \cite{t:cm}.

\ble\label{lem:com}
The $(i,j)$ entry of the matrix ${\wh L}^{\pm}(u)$ is given by
\beql{mateco}
{\wh l}^{\ts\pm}_{ij}(u)=(-q)^{j-i}\ts
{l^{\pm}\ts}^{1\ts\dots\ts
\wh j\ts\dots\ts n}_{1\ts\dots\ts\wh i\ts\dots\ts n}(u),
\eeq
where the hats on the right-hand side
indicate the indices to be omitted.
\ele

\bpf
According to \eqref{qminorgen},
the quantum determinants are defined by the relations
\ben
A^{(n)}_q L^{\pm}_1(uq^{2n-2})\cdots L^{\pm}_n(u)=A^{(n)}_q\tss \qdet L^{\pm}(u).
\een
Hence, the definition \eqref{comatrix} of the quantum comatrices implies
\ben
A^{(n)}_q L^{\pm}_1(uq^{2n-4})\cdots L^{\pm}_{n-1}(u)=A^{(n)}_q\tss {\wh L}^{\pm}_n(u).
\een
Apply both sides to the vector $e_1\ot\dots\ot \wh e_i\ot\dots\ot e_n\ot e_j$
and use \eqref{pq} and \eqref{antisym} to equate the coefficients of
the vector $A^{(n)}_q(e_1\ot\dots\ot e_n)$
to obtain \eqref{mateco}.
\epf

We will also use the relations for the transposed comatrices; see \cite{nt:yg} and \cite{t:cm}.

\ble\label{lem:transco}
We have the relations
\ben
D{\wh L}^{\pm}(u)^tD^{-1}L^{\pm}(uq^{2n-2})^t=\qdet L^{\pm}(u)\ts 1.
\een
\ele

\bpf
Since each relation only involves generators belonging to the
subalgebra of $\U_q(\wh\gl_n)$ generated
by the coefficients of the series $l_{ij}^+(u)$ or $l_{ij}^-(u)$,
we may assume that the central charge is specialized to zero, $c=0$,
and derive the desired relations in the quotient algebra $\U^{\circ}_q(\wh\gl_n)$.
The mapping
\ben
\theta: L^{\pm}(u)\to L^{\mp}(u^{-1})^t
\een
defines an automorphism of the algebra $\U^{\circ}_q(\wh\gl_n)$.
By Lemma~\ref{lem:com} and the quantum minor formulas,
for the images under $\theta$ we have
\ben
\theta:{\wh L}^{\pm}(u)\mapsto D{\wh L}^{\mp}(u^{-1}q^{-2n+4})^tD^{-1}
\een
and
\ben
\theta: \qdet L^{\pm}(u)\mapsto \qdet L^{\mp}(u^{-1}q^{-2n+2}).
\een
It remains to apply $\theta$ to both sides of relations \eqref{comatrix}
and then replace $u$ by $u^{-1}q^{-2n+2}$.
\epf

Now, using \eqref{comatrix} and the centrality of $\qdet L^{\pm}(u)$, we can write
\eqref{lltr} as
\ben
z^{\pm}(u)=L^{\pm}(uq^{2n})^tD(L^{\pm}(u)^{-1})^tD^{-1}
=(\qdet L^{\pm}(u))^{-1}L^{\pm}(uq^{2n})^tD{\wh L}^{\pm}(uq^2)^tD^{-1}.
\een
By applying Lemma~\ref{lem:transco}
with $u$ replaced by $uq^2$ we find that this expression coincides
with the right hand side of \eqref{liouvpm}.
\epf

\bre\label{rem:central}
Since the second part of Proposition~\ref{prop:equiz} was not used in the proof of Theorem~\ref{thm:liov},
the fact that the coefficients of both series $z^+(u)$ and $z^-(u)$ belong to the center of
the quantum affine algebra $\U_q(\wh\gl_n)$ also follows from
\eqref{liouvpm} due to the respective properties
of the quantum determinants $\qdet L^{+}(u)$ and $\qdet L^{-}(u)$.
\ere

\section{Quantum Gelfand invariants}
\label{sec:qgi}

We will follow \cite{df:it} to define the {\em quantized enveloping algebra\/} $\U_q(\gl_n)$ in its
$R$-matrix presentation \cite{rtf:ql} as the algebra
generated by elements $l^{+}_{ij}$ and $l^{-}_{ij}$ with $1\leqslant i,j\leqslant n$
subject to the relations
\ben
\bal
l^-_{ij}&=l^+_{ji}=0, \qquad 1 \leqslant i<j\leqslant n,\\
l^-_{ii}\ts l^+_{ii}&=l^+_{ii}\ts l^-_{ii}=1,\qquad 1\leqslant i\leqslant n,\\[0.3em]
R\ts L^{\pm}_1L^{\pm}_2&=L^{\pm}_2L^{\pm}_1R,\qquad
R\ts L^+_1L^-_2=L^-_2 L^+_1R,
\eal
\een
where the $R$-matrix $R$ is defined in \eqref{R}, while
$L^+$ and $L^-$ are the matrices
\ben
L^{\pm}=\sum_{i,j}e_{ij}\ot l^{\pm}_{ij}\in\End\CC^n\ot \U_q(\gl_n).
\een
The universal enveloping algebra $\U(\gl_n)$
is recovered from $\U_q(\gl_n)$ in the limit
$q\to 1$ by the formulas
\begin{alignat}{2}\label{taugen1}
\frac{l^-_{ij}}{q-\qin}&\rightarrow E_{ij},
\qquad\quad
\frac{l^+_{ji}}{q-\qin}&&\rightarrow -E_{ji}\qquad\quad\text{for}\quad i> j,\\
\intertext{and}
\label{taugen2}
\frac{l^-_{ii}-1}{q-1}&\rightarrow E_{ii},
\qquad\quad
\frac{l^+_{ii}-1}{q-1}&&\rightarrow -E_{ii}\qquad\quad\text{for}\quad i=1,\dots,n.
\end{alignat}

\medskip

Set $M=L^-(L^+)^{-1}$.
By \cite{rtf:ql}, the quantum traces defined by
$
\tr_q\tss M^m=\tr\ts DM^m
$
belong to the center of the algebra $\U_q(\gl_n)$.
The elements $\tr_q\tss M^m$ act by multiplication
by scalars in the irreducible highest weight representations $L_q(\la)$.
The representation
$L_q(\la)$ of $\U_q(\gl_n)$
is generated by a nonzero vector $\xi$ such that
\begin{alignat}{2}
l^+_{ij}\ts\xi&=0 \quad &&\text{for} \quad
1\leqslant i<j\leqslant n,
\non\\
l^-_{ii}\ts\xi&=q^{\la_i}\tss\xi \quad
&&\text{for} \quad 1\leqslant i\leqslant n,
\non
\end{alignat}
for an $n$-tuple $\la=(\la_1,\dots,\la_n)$ of integers (or real numbers).
This is a $q$-deformation of the irreducible
$\gl_n$-module $L(\la)$
with the highest weight $\la$.

\paragraph{Proof of Theorem~\ref{thm:eigen}.} We will show that
the eigenvalue of the quantum Gelfand invariant $\tr_q\tss M^m$ in $L_q(\la)$ is given by
formula \eqref{qppexp}.

Recall the evaluation homomorphism $\U_q(\wh\gl_n)\to \U_q(\gl_n)$
defined by
\beql{eval}
L^+(u)\mapsto L^+-L^- u,\qquad
L^-(u)\mapsto L^--L^+ u^{-1},\qquad q^c\mapsto 1,
\eeq
and apply it to both sides of the Liouville formula
\beql{liuplus}
z^{+}(u)=\frac{\qdet L^{+}(uq^{2})}{\qdet L^{+}(u)}
\eeq
proved in Theorem~\ref{thm:liov}. The image of the quantum
determinant $\qdet L^{+}(u)$ is found by
\ben
\sum_{\si\in \Sym_n} (-q)^{-l(\si)} \ts \big(l^{+}_{\si(1)1}-l^-_{\si(1)1}uq^{2n-2}\big)\cdots
\big(l^{+}_{\si(n)n}-l^-_{\si(n)n}u\big).
\een
By applying this element to the highest vector $\xi$ of $L_q(\la)$, we find that
its eigenvalue is given by
\ben
\big(q^{-\la_1}-q^{\la_1+2n-2}u\big)\dots \big(q^{-\la_n}-q^{\la_n}u\big)
=q^{n(n-1)/2}(q^{-\ell_1}-q^{\ell_1}u)\dots (q^{-\ell_n}-q^{\ell_n}u)
\een
with $\ell_i=\la_i+n-i$. Therefore, the eigenvalue of the image
of the right hand side of \eqref{liuplus} in $L_q(\la)$ is given by
\ben
\frac{(q^{-\ell_1}-q^{\ell_1+2}u)\dots (q^{-\ell_n}-q^{\ell_n+2}u)}
{(q^{-\ell_1}-q^{\ell_1}u)\dots (q^{-\ell_n}-q^{\ell_n}u)}.
\een
To expand this rational function into a series in $u$ write it as
\ben
C+\frac{a_1}{1-q^{2\ell_1}u}+\dots+\frac{a_n}{1-q^{2\ell_n}u}
\een
to find that the constants $a_k$ are given by
\ben
a_k=(q^{n-1}-q^{n+1})\ts \frac{[\ell_1-\ell_k+1]_q\dots [\ell_n-\ell_k+1]_q}
{[\ell_1-\ell_k]_q\ldots\wedge\dots [\ell_n-\ell_k]_q}.
\een
Then write
\ben
\frac{1}{1-q^{2\ell_k}u}=\sum_{m=0}^{\infty}\ts q^{2\ell_k m}\tss u^m.
\een

On the other hand, using \eqref{zu}, for the image of $z^+(u)$
under the evaluation homomorphism \eqref{eval} we get
\ben
\frac{1}{[n]_q}\ts \tr\ts D(L^{+}-L^-uq^{2n})(L^+-L^-u)^{-1}
=\frac{1}{[n]_q}\ts \tr\ts D(1-Muq^{2n})(1-Mu)^{-1}
\een
which equals
\ben
1+(q^{n-1}-q^{n+1})\ts\sum_{m=1}^{\infty}\ts \tr_q\tss M^m\tss u^m.
\een
Formula \eqref{qppexp} now follows by equating the coefficients of the powers $u^m$
on both sides of the power expansions. Note that the formula is also valid for $m=0$.
\qed

\bre\label{rem:pp}
Due to \eqref{cllim}, the Gelfand invariant $\tr\ts E^m\in\U(\gl_n)$ is
obtained as the limit value as $q\to 1$ of the expression
\ben
\frac{1}{(q-\qin)^m}\ts\tr\ts D(M-1)^m=
\frac{1}{(q-\qin)^m}\ts\sum_{r=0}^m\binom{m}{r}(-1)^{m-r}\ts\tr_q\tss M^r.
\een
We will look at the limit value of the corresponding linear combination of the eigenvalues
given by \eqref{qppexp} and note that
the $q$-numbers specialize by the rule $[r]_q\to r$ as $q\to 1$.
Hence, for each $k=1,\dots,n$ it suffices to
find the limit value of the expression
\ben
\frac{1}{(q-\qin)^m}\ts\sum_{r=0}^m\binom{m}{r}(-1)^{m-r}\ts q^{2\ell_k r}
=\frac{(q^{2\ell_k}-1)^m}{(q-\qin)^m},
\een
which equals $\ell_k^{\tss m}$, thus yielding formula \eqref{ppexp}.
\qed
\ere

\medskip

Since the evaluation homomorphism \eqref{eval} is surjective, the centrality property
of the elements $\tr_q\tss M^m$ also follows from the above calculations due to
Proposition~\ref{prop:equiz}. Similarly, by using the three remaining formulas
in \eqref{zu} and \eqref{opzu}, we get three more families of central
elements in $\U_q(\gl_n)$ together with the relations between them given by
\ben
\tr\ts D^{-1}((L^+)^{-1}L^-)^m=\tr_q\tss M^m\Fand
\tr\ts D^{-1}((L^-)^{-1}L^+)^m=\tr\ts D(L^+(L^-)^{-1})^m.
\een
Moreover, additional relations between the families are provided by
the isomorphism
\ben
\U_q(\gl_n)\to\U_{\qin}(\gl_n),\qquad L^{\pm}\mapsto (L^{\pm})^{-1}.
\een
This implies that the formulas for the eigenvalues of the other two central elements
of $\U_q(\gl_n)$ in $L_q(\la)$ are obtained by the replacement $q\mapsto\qin$
in \eqref{qppexp}:
\ben
\tr\ts D(L^+(L^-)^{-1})^m\mapsto\sum_{k=1}^n\ts q^{-2 \ell_k m}\ts
\frac{[\ell_1-\ell_k+1]_q\dots [\ell_n-\ell_k+1]_q}
{[\ell_1-\ell_k]_q\ldots\wedge\dots [\ell_n-\ell_k]_q}.
\een
Alternatively, these eigenvalue formulas can be derived in the same way as in the proof
of Theorem~\ref{thm:eigen} by working with the other Liouville formula in \eqref{liouvpm}
instead of \eqref{liuplus}.

\section*{Declarations}

\subsection*{Competing interests}
The authors have no competing interests to declare that are relevant to the content of this article.

\subsection*{Funding}
This work was supported by the Australian Research Council, grant DP240101572.

\subsection*{Availability of data and materials}
No data was used for the research described in the article.

\bigskip\bigskip

\small

\noindent
N.J.:\qquad\qquad\qquad\qquad\\
Department of Mathematics\\
North Carolina State University, Raleigh, NC 27695, USA\\
jing@math.ncsu.edu

\vspace{3 mm}

\noindent
M.L.:\newline
School of Mathematical Sciences\\
South China University of Technology\\
Guangzhou, Guangdong 510640, China\\
mamliu@scut.edu.cn

\vspace{3 mm}

\noindent
A.M.:\newline
School of Mathematics and Statistics\newline
University of Sydney,
NSW 2006, Australia\newline
alexander.molev@sydney.edu.au

\end{document}